\documentclass[11pt,leqno]{article}

\usepackage{amstex,amscd,theorem}

\long\def\comment#1\endcomment{}

\makeatletter
\begingroup
\gdef\th@dotted{\normalfont\itshape
  \def\@begintheorem##1##2{%
        \item[\hskip\labelsep \theorem@headerfont ##1\ ##2.]}%
\def\@opargbegintheorem##1##2##3{%
   \item[\hskip\labelsep \theorem@headerfont ##1\ ##2\ (##3).]}}
\endgroup
\makeatother

\theoremstyle{dotted}

\newtheorem{theorem}{Theorem}[section]
\newtheorem{lemma}[theorem]{Lemma}
\newtheorem{corr}[theorem]{Corollary}
\newtheorem{prop}[theorem]{Proposition}

\newtheorem{claim}[theorem]{Claim}

\makeatletter
\begingroup
\gdef\th@upshape{\normalfont
  \def\@begintheorem##1##2{%
        \item[\hskip\labelsep \theorem@headerfont ##1\ ##2.]}%
\def\@opargbegintheorem##1##2##3{%
   \item[\hskip\labelsep \theorem@headerfont ##1\ ##2\ (##3).]}}
\endgroup
\makeatother

\theoremstyle{upshape}

\newtheorem{defn}[theorem]{Definition}
\newtheorem{remark}[theorem]{Remark}

\makeatletter
\@addtoreset{equation}{section}
\makeatother

\newcommand{\proof}[1][Proof.]{\smallskip\noindent{\em #1}}
\def\endproof{\hfill\ensuremath{\square}\par\medskip}

\def\eqref#1{\thetag{\ref{#1}}}

\let\latexref=\ref
\def\ref#1{{\normalfont{\latexref{#1}}}}

\newcommand{\wt}{\widetilde}
\newcommand{\wh}{\widehat}

\newcommand{\C}{{\Bbb C}}
\newcommand{\R}{{\Bbb R}}
\newcommand{\Z}{{\Bbb Z}}
\newcommand{\Q}{{\Bbb Q}}

\newcommand{\A}{{\Bbb A}}

\newcommand{\cp}{{\Bbb CP}}

\newcommand{\calo}{{\cal O}}
\newcommand{\E}{{\cal E}}
\newcommand{\T}{{\cal T}}

\newcommand{\X}{{\frak X}}
\newcommand{\Y}{{\cal Y}}
\newcommand{\m}{{\frak m}}
\newcommand{\TT}{{\frak t}}
\newcommand{\hh}{{\frak h}}

\newcommand{\SL}{{\frak s}{\frak l}}

\newcommand{\Spec}{\operatorname{Spec}}
\newcommand{\Proj}{\operatorname{Proj}}
\newcommand{\Bl}{\operatorname{Bl}}           

\newcommand{\Pic}{\operatorname{Pic}}

\newcommand{\Aut}{\operatorname{Aut}}

\newcommand{\id}{\operatorname{\sf id}} 
\newcommand{\rk}{\operatorname{\sf rk}} 
\newcommand{\cchar}{\operatorname{\sf char}} 
\renewcommand{\dim}{\operatorname{\sf dim}} 
\newcommand{\codim}{\operatorname{\sf codim}} 
\newcommand{\cl}{\operatorname{\sf cl}}
\newcommand{\res}{\operatorname{\sf res}}

\title{Dynkin diagrams and crepant resolutions of quotient singularities}

\author{D. Kaledin}

\date{}

\begin{document}

\maketitle

\tableofcontents

\section*{Introduction}

Let $V$ be a complex vector space, and assume that a finite group
$G$ acts on $V$ by linear transformations. Consider the quotient
$V/G$. If $\dim V = 1$, the quotient $V/G$ is a smooth algebraic
variety. If $\dim V \geq 2$, it may still happen that the quotient
$V/G$ is smooth. However, usually the algebraic variety $V/G$ is
singular. In this case it is natural to look for good resolutions of
singularities $X \to V/G$ of the quotient $V/G$.

There exists a well-known situation when there indeed exists such a
resolution: the so-called {\em McKay correspondence}. This is the
case when the vector space $V = \C^2$ is of dimension $\dim V = 2$,
and the finite subgroup $G \subset SL(2,\C)$ preserves the standard
symplectic form on $V$. For such a quotient $V/G$, one can construct
a canonical smooth resolution $X \to V/G$ which has many good
properties. One of these properties is that the smooth algebraic
variety $X$ has trivial canonical bundle. Such resolutions are
called {\em crepant}.

In their recent paper \cite{hilbres}, Y. Ito and H. Nakajima studied
possible generalizations of the McKay correspondence to higher
dimensions. One of the questions that they asked was the following:
for which pairs $\langle V, G\rangle$ the quotient $V/G$ admits a
crepant smooth resolution $X \to V/G$?

Aside from the McKay correspondence, the best-known example of such
a pair is the $2n$-dimensional complex vector space $V = \C^{2n}$
and the symmmetric group $G=S_n$ on $n$ letters that acts on $V$ by
transpositions, separately on odd and even coordinates. In this
case, a crepant resolution $X$ is the Hilbert scheme of $n$ points
on $\C^2$ (see \cite{Nhilb}).

We note that both in the case of the McKay correspondence and in the
case of the Hilbert scheme, the complex vector space $V$ carries a
symplectic form preserved by the group $G$. This suggests that the
natural first choice for a quotient that admits a crepant resolution
would be the quotient $V/G$ of a vector space $V$ which carries a
$G$-invariant symplectic form.

Such symplectic vector spaces are, of course, very common. A simple
way to construct them is by ``doubling'' representations of finite
groups: for an arbitrary vector space $V_o$ equipped with a linear
action of a finite group $G$, the sum $V = V_o \oplus V_o^*$ of the
vector space $V_o$ with its dual $V_o^*$ carries both a natural
$G$-action and a canonical $G$-invariant symplectic form.

\medskip

In this paper, we study smooth crepant resolutions of quotients
$V/G$, where $V = V_o \oplus V_o^*$ is a symplectic vector space
obatined by doubling a representation $V_o$ of a finite group $G$.
The main result of this paper, Theorem~\ref{main}, claims that if
for $V = V_o \oplus V_o^*$ the quotient $V/G$ admits a smooth
crepant resolution, then the action of the group $G$ on the vector
space $V_o$ is generated by (complex) reflections.

\medskip

Finite groups $G \subset \Aut V_o$ of automorphisms of a complex
vector space $V_o$ which are generated by reflections have been an
object of much study, and there exists a complete classification of
pairs $\langle V_o, G \rangle$ of this type. The best-known part of
this classification concerns subgroups $G \subset \Aut V_o$ which
preserve a rational lattice $V_\Q \subset V_o$. In this case $G$ is
generated by reflections if and only if it is a product of Weyl
groups associated to Dynkin diagrams of finite type (see
\cite[Ch. VI, \S 2, p. 5, Proposition 9]{Bourbaki}).

Of all the subgroups $G \subset \Aut(\C^2)$ studied in the McKay
correspondence, only the simplest ones preserve the decompostion
$\C^2 = \C \oplus \C$, thus falling under the assumptions of our
main theorem. These are the cyclic groups $G = \Z/n\Z$. From the
point of view of the theory of groups generated by reflections, this
is the trivial case. The case of the Hilbert scheme is more
interesting. It also falls under our assumptions, and corresponds to
the Weyl group of the type $A_n$.

Among the quotients by the Weyl groups of other types, a crepant
resolution has been constructed for the type $C_n$ by A. Kuznecov
(\cite{K}, see also \cite{part-res}). For the rest of the Dynkin
diagrams the question posed in \cite{hilbres} seems to be open.

Aside from our main theorem, we also obtain some results on the
structure of smooth crepant resolutions of quotient sungularities
$V/G$, and we prove some additional facts. Among these, we would
like to mention Theorem~\ref{equivar}. It claims, more or less, that
if a connected algebraic group $H$ acts on the space $V$, and if the
$H$-action commutes with the $G$-action, then for any smooth crepant
resolution $\pi:\wt{X} \to X$ the induced $H$-action on the quotient
$X = V/G$ lifts to an $H$-action on the variety $\wt{X}$. We refer
the reader to Section 1 for the precise statement.

For this result, we do not need the assumption $V = V_o \oplus
V_o^*$. It holds for an arbitrary symplectic vector space $V$ and a
finite subgroup $G \subset \Aut V$ which preserves the symplectic
form. In fact, even the symplectic form is not necessary -- it
suffices to assume that the group $G$ preserves a volume form on the
space $V$ (see Remark~\ref{CalabiYau}).

Another result which we would like to mention here is the uniqueness
statement Theorem~\ref{unique}. It claims that if any two
reflections in the subgroup $G \subset \Aut V_o$ are conjugate
within $G$, then there exists at most one smooth crepant resolution
of the quotient $V/G$. When $G \subset \Aut V_o$ is a Weyl group,
this condition is satisfied if and only if the corresponding Dynkin
diagram is simply laced.

Unfortunately, for technical reasons we need to impose an additional
restriction: in Theorem~\ref{unique} we require that the complex
vector space $V_o$ admits a real structures invariant under
$G$. Among other things, this implies that $V_o \cong V_o^*$ as
representations of the group $G$.

\medskip

\noindent
Let us now give a brief outline of the contents of the paper.
\begin{itemize}
\item Section 1 contains the precise description of the general
setup, the definitions, and the statements of all the results proved
in the paper.

\item A very short Section 2 gives an overview of the simplest
possible case $V = \C^2$, $G = \{\pm 1\}$.

\item Section 3 is devoted to some general technical results -- in
particular, we prove that under our assumptions, for every smooth
crepant resolution $X \to V/G$ the variety $X$ is holomorphically
symplectic. The reader is advised to skip this section at first
reading.

\item In Section 4, we introduce a certain stratification of the
quotient $V/G$.

\item This stratification is then used in Section 5 to study in some
detail the geometry of an arbitrary smooth crepant resolution $X \to
V/G$. In particular, we compute the rational Picard group $\Pic(X)
\otimes \Q$. After that we are able to prove most of our results --
in fact all of them, except for the main Theorem~\ref{main}.

\item In Section 6 we introduce, following \cite{quiver}, a certain
canonical action of the group $\C^*$ on an arbitrary smooth crepant
resolution $X \to V/G$. Using this action, we then prove
Theorem~\ref{main}.
\end{itemize}

To finish the introduction, I would like mention that while this
paper was in preparation, a stronger version of Theorem~\ref{main}
was proved by M. Verbitsky. This will be the subject of his upcoming
paper \cite{Vnew}. His proof is different from the one presented
here; it is much simpler and it does not use most of our technical
results. Thus the last section of this paper is perhaps obsolete. I
have decided to keep it anyway, since the $\C^*$-action which is
constructed and studied there may be of independent interest.

\medskip

\noindent
{\bf Acknowledgements.} I would like to thank R. Bezrukavnikov,
F. Bogomolov, V. Ginzburg, A. Kuznecov, M. Leenson and M. Reid for
very helpful discussions, either in person or by e-mail. I am
particularly grateful to Misha Verbitsky and Tony Pantev for many
valuable suggestions and important insights.

\section{Preliminaries}\label{prelim}

Let $V$ be a complex vector space of dimension $2n$ equipped with a
non-degenerate symplectic form $\Omega_V$.  Assume given a finite
group $G \subset \Aut V$ of automorphisms of the vector space $V$
which preserves the symplectic form.  The quotient $V/G$ is then
naturally a singular affine algebraic variety over $\C$.  In this
paper we will interested in smooth resolutions of the variety $V/G$.

To fix terminology, we introduce the following.

\begin{defn}
A map $\pi:\wt{X} \to X$ is called a {\em resolution} of an
irreducible algebraic variety $X$ if the algebraic variety $\wt{X}$
is irreducible, and the map $\pi$ is generically one-to-one. A
resolution $\pi:\wt{X} \to X$ is called {\em smooth} if the
algebraic variety $\wt{X}$ is smooth, and {\em projective} if the
map $\pi:\wt{X} \to X$ is projective.
\end{defn}

We will only be interested in a special class of resolutions, the
so-called {\em crepant} ones. The general definition is as follows.
Let $\pi:\wt{X} \to X$ be a resolution of an irreducible algebraic
variety $X$. Assume that both $X$ and $\wt{X}$ are normal and admit
canonical bundles $K_X$ and $K_{\wt{X}}$\footnote{An algebraic
variety $X$ is said to admit a canonical bundle if the canonical
bundle of the non-singular part of $X$ extends to a line bundle on
the whole $X$. If $X$ is normal, then such an extension is
unique}. Let $U \subset X$ be a non-singular open dense subset such
that $\pi:\wt{U} = \pi^{-1}(U) \to U$ is one-to-one.

\begin{defn}
The resolution $\pi:\wt{X} \to X$ is called {\em crepant} if the
canonical isomorphism
$$
\pi^*K_U \cong K_{\wt{U}},
$$
over $U = \pi^{-1}(U) \subset \wt{X}$ extends to a bundle
isomorphism $\pi^*K_X \cong K_{\wt{X}}$ over the whole $\wt{X}$.
\end{defn}

The quotient variety $X = V/G$ is normal and irreducible. Moreover,
it obviously admits a canonical bundle $K_X$ -- namely, the trivial
one. Therefore for every crepant resolution $\pi:\wt{X} \to X$ the
canonical bundle $K_{\wt{X}}$ is also trivial. 

For projective resolutions, this necessary condition is in fact
sufficient: a projective resolution $\pi:\wt{X} \to X$ is crepant
if the canonical bundle $K_{\wt{X}}$ is trivial.  Indeed, in this
case a map
$$
\pi^*K_X \cong \calo_{\wt{X}} \to K_{\wt{X}} \cong \calo_{\wt{X}}
$$
is simply a function on $\wt{X}$. We are given a function $f$
without zeroes on $\wt{U} = \pi^{-1}(U) \subset \wt{X}$, and we have
to prove that it extends to a function without zeroes on the whole
$\wt{X}$. It suffices to extend $\pi^{-1} \circ f:U \to \C$ to a
function $f$ without zeroes on $X$ and take $\pi^*f$. But since $X$
is normal and the map $\pi:\wt{X} \to X$ is dominant, we can assume
that the complement $X \setminus U \subset X$ to the open
non-singular subset $U \subset X$ is of codimension $\codim X
\setminus U \geq 2$. Again using the fact that $X$ is normal, we see
that every function $f$ without zeroes on $U$ extends to a function
without zeroes on the whole $X$.

In this paper we consider a quotient $X=V/G$ of a symplectic vector
space $V$ which admits a smooth projective crepant resolution
$\pi:\wt{X} \to X$, and we prove assorted results on the geometry of
the quotient $X$ which hold under varying linear-algebraic
assumptions on the pair $\langle V,G \rangle$. Our first result
holds for every symplectic vector space.

\begin{theorem}\label{equivar}
Let $X = V/G$ be the quotient of a symplectic vector space $V$ by a
symplectic action of a finite group $G \subset \Aut(V)$. Assume
given a smooth projective crepant resolution $\pi:\wt{X} \to X$.
\begin{enumerate}
\item For every $G$-invariant vector field $\TT_V$ on $V$, the
induced vector field $\TT_0$ on the non-singular part $X_0 \subset X
= V/G$ lifts to a vector field $\TT$ on the whole smooth algebraic
variety $\wt{X}$.
\item If a connected algebraic group $H$ acts on the algebraic
variety $V$, and if the $H$-action commutes with the $G$-action,
then the induced $H$-action on the quotient $X = V/G$ lifts to an
$H$-action on the variety $\wt{X}$.
\item The same is true for an arbitrary, not necessarily connected
algebraic group $H$, which acts on the vector space $V$ by linear
transformations.
\end{enumerate}
\end{theorem}

Our second general result is somewhat technical, but it might be
useful in applications. To formulate this result, consider the
formal completion $\wh{X}$ of a quotient variety $X = V/G$ at $0 \in
V/G$. Say that a formal scheme $\X$ equipped with a projective map
$\pi:\X \to \wh{X}$ is a smooth crepant resolution of the completion
$\wh{X}$ if $\X$ is smooth and the canonical bundle $K_{\X}$ is
trivial.

\begin{theorem}\label{formal}
Assume given a smooth projective crepant resolution $\pi:\X \to
\wh{X}$ of the completion $\wh{X}$. Then there exists a unique
smooth projective crepant resolution $\pi:\wt{X} \to X$ which gives
$\pi:\X \to \wh{X}$ after completing along the fiber $\pi^{-1}(0)
\subset \wt{X}$.
\end{theorem}

\begin{remark}
By virtue of Grothendieck's algebraization theorem
\cite[Th\'eor\`eme 5.4.5]{EGA3}, one can replace here a formal
scheme $\X$ with a usual scheme projective over $\wh{X}$.
\end{remark}

Our final result which holds for every symplectic vector space $V$
is the following.

\begin{theorem}\label{inductive}
Assume that the quotient $X = V/G$ of a symplectic vector space $V$
by a symplectic action of a finite group $G \subset \Aut(V)$ admits
a smooth projective crepant resolution. Let $v \in V$ be an
arbitrary vector, let $G_v \subset G$ be the stabilizer of the
vector $v \in V$, and let $V' \subset V$ be the unique
$G_v$-invariant complement to the subspace $V^{G_v} \subset V$ of
vectors fixed by $G_v$. 

Then the quotient $X_v = V'/G_v$ also admits a smooth projective
crepant resolution.
\end{theorem}

The main result of the paper holds only under an additional
assumption. To formulate this result, recall that a finite-order
automorphism $g:V \to V$ of a vector space $V$ is called a {\em
complex reflection} if the subspace $V^g \subset V$ of
$g$-invariants has codimension exactly $1$.

\begin{theorem}\label{main}
Let $V$ be a symplectic vector space, and let $G \subset \Aut(V)$ be
a finite group of symplectic automorphims of the vector space
$V$. Assume that the symplectic vector space $V$ admits a
$G$-invariant Lagrangian subspace $V_o \subset V$. 

If the quotient $X = V/G$ admits a smooth projective crepant
resolution $\pi:\wt{X} \to X$, then the subgroup $G \subset
\Aut(V_o)$ is generated by complex reflections.
\end{theorem}

Note that for an arbitrary representation $V_o$ of a finite group
$G$, the vector space $V = V_o \oplus V_o^*$ is symplectic and
satisfies the conditions of Theorem~\ref{main}.

Theorem~\ref{main} is analogous to the following classic fact.

\def\attrib{\cite[\!Ch.\ \!\!V, \S 5, Theorem 4]{Bourbaki}}
\begin{theorem}[\attrib]\label{bourb}
The quotient $V/G$ \!of a complex vector space $V$ by a finite
subgroup $G \subset \Aut(V)$ is smooth if and only if the subgroup $G
\subset \Aut(V)$ is generated by complex reflections. \endproof
\end{theorem}

In fact, our proof of Theorem~\ref{main} uses this fact
directly. Note, however, that we do not claim the converse to
Theorem~\ref{main}. At present, we do not know for which of the
subgroups $G \subset \Aut(V_o)$ generated by complex reflections the
quotient $V/G$ admits a smooth crepant resolution. We hope to return
to this in a later paper.

In the course of proving Theorem~\ref{main}, we also establish the
following uniqueness statement, which might be of independent
interest. Unfortunately, we can only prove this result under an even
stronger additional assumption.

\begin{theorem}\label{unique}
In the notation of Theorem~\ref{main}, assume that the complex
vector space $V_o$ admits a real structure preserved by the group
$G$. Assume also that the reflections in $G$ form a single conjugacy
class.

Then every two smooth projective crepant resolutions $\pi_1:\wt{X}_1
\to X$, $\pi_2:\wt{X}_2 \to X$ of the quotient $V/G$ are canonically
isomorphic.
\end{theorem}

In fact, the smoothness assumption on $\wt{X}_1$, $\wt{X}_2$ imposed
in this theorem can be weakened (see Remark~\ref{semism.unique}).

\section{The case of dimension $1$}

We begin with the simple and well-known case $\dim V = 2$. In this
case Theorem~\ref{inductive} and Theorem~\ref{main} are trivial. On
the other hand, the proofs of Theorem~\ref{equivar} and
Theorem~\ref{formal} are essentially the same as in the general
case. We postpone these till Section~\ref{sec.pic}. In this section
we concentrate on the uniqueness result, Theorem~\ref{unique}.

Under the additional assumptions of Theorem~\ref{unique}, the group
$G$ preserves a $1$-dimensional Lagrangian subspace $V_o \subset V$
and a real structure on $V_o$. Therefore $G$ necessarily consists of
multiplications by $\pm 1$. The affine quotient variety $X = V
/\{\pm 1\}$ is the singular quadric in $\C[u,v,w]$ given by the
equation $uw=v^2$, and the quotient map $V = \C[x,y] \to X$ sends
$u$ to $x^2$, $v$ to $xy$ and $w$ to $y^2$.

A classic crepant resolution for $X$ is obtained by blowing up $0
\subset X$. It coincides with the total space of the cotangent
bundle to the complex projective line $\cp^1$.

We will now prove Theorem~\ref{unique} for this simple case. Assume
given another smooth projective crepant resolution $\pi:\wt{X} \to
X$.  By definition the map $\pi$ is one-to-one over a non-singular
open subset $X_0 \subset X$, so that $\pi:\wt{X}_0 = \pi^{-1}(X_0)
\to X_0$ is an isomorphism. Since $X = V/G$ is normal, we can assume
that the complement $X \setminus X_0$ is of codimension $\codim X
\setminus X_0 > 1$. Moreover, we can assume that $X_0$ does not
contain $0 \in X$, so that the quotient map $V \to X$ is \'etale
over $X_0 \subset X$.

The $G$-invariant symplectic form $\Omega_V$ on $V$ induces a
non-degenerate symplectic form $\Omega$ on $X_0 \subset X$. But
$\dim X = 2$. Therefore a symplectic form on $X_0$ is a section of
the canonical bundle $K_X$. Since the resolution $\pi:\wt{X} \to X$
is crepant, the pullback $\pi^*\Omega$ extends to a non-degenerate
symplectic form $\Omega$ on the whole smooth algebraic variety
$\wt{X}$. The $2$-form $\Omega$ on $\wt{X}$ induces a canonical
isomorphism
$$
\omega:H^0(\wt{X},\T) \overset{\sim}{\to} H^0(\wt{X},\Omega^1)
$$
between the space of global $1$-forms on $\wt{X}$ and the space of
global vector fields on $\wt{X}$.

The standard action of the group $SL(2,\C)$ on the vector space $V =
\C^2$ commutes with the action of the group $G$. Therefore it
descends to an action on the space $X$. The corresponding
infinitesemal action is generated by three vector fields $E,F,H$ on
$X_0 \cong \wt{X}_0$ which generate the Lie algebra $\SL(2)$. Over
$X_0$, the isomorphism $\omega$ sends these three vector fields to
the differential forms $du$, $dv$ and $dw$.

Since these forms obviously extend to the whole $\wt{X}$, we
conclude that the vector fields $E$, $F$ and $H$ also extend to the
whole $\wt{X}$.

By assumption the map $\pi:\wt{X} \to X$ is projective. Therefore
the variety $\wt{X}$ carries a line bundle $L$ which is very ample
for the map $\pi$, so that we have
$$
\wt{X} = \Proj\bigoplus_{k \geq 0}\pi_*L^k.
$$
Let $j:X_0 \hookrightarrow X$ be the embedding map of the open
subset $X_0 = X \setminus \{0\} \subset X$. Since the codimension
$\codim X \setminus X_0$ is greater than $1$, the Picard group of
the open subset $\wt{X}_0 \cong X_0 \subset X$ coincides with the
Picard group of the whole variety $X$, which is isomorphic to
$\Z/2\Z$. Replacing $L$ with its positive power, we can assume that
$L|_{\wt{X}_0}$ is trivial.  Any isomorphism $\pi_*L
\overset{\sim}{\rightarrow} \calo_X$ over $X_0$ extends to an
embedding
$$
\pi_*L \hookrightarrow j_*\calo_X \cong \calo_X.
$$
This means that the sheaf $\pi_*L$ is in fact a sheaf of ideals in
$\calo_X$, and the resolution $\pi:\wt{X} \to X$ is the blowup of
the variety $X$ in this sheaf of ideals.

Since the $\SL(2)$-action on $X$ lifts to an action on $\wt{X}$, the
sheaf of ideals $\pi_*L \subset \calo_X$ must be
$\SL(2)$-invariant. Equivalently, the ideal of global sections
$H^0(X,\pi_*L) \subset H^0(\calo_X)$ must be $\SL(2)$ invariant. But
the only $\SL(2)$-invariant ideals in the algebra $H^0(\calo_X) =
\C[u,v,w]/(uw-v^2)$ are the powers of the maximal ideal
corresponding to the point $0 \in X$. Consequently, $\wt{X}$ must be
isomorphic to the blow-up of the point $0 \in X$.

This argument proves Theorem~\ref{unique} in the case $V = \C^2$.
We sum up the results for this simple situation in the following
statement.

\begin{prop}\label{dim.1}
Every projective smooth crepant resolution $\pi:\wt{X} \to X$ is
isomorphic to the blow-up of the point $0 \in X$. Moreover, the
Picard group of line bundles on $\wt{X}$ is $\Z$, and for the ample
line bundle corresponding to $2k \in \Z$, $k \geq 1$ we have a
canonical isomorphism
$$
\pi_*L \cong \m^k,
$$
where $\m \subset \calo_X$ is the sheaf of ideals of the point $0
\in X$. \endproof
\end{prop}

Note that everything in the proof carries over literally to the
formal situation of Theorem~\ref{formal}. Therefore the same
statements hold for projective smooth crepant resolutions of the
formal completion $\wh{X}$ of $X$ at $0 \subset X$.

\section{The holomorphic symplectic form}

In higher dimensions the canonical bundle is no longer the bundle of
$2$-forms. Still, it turns out that for every symplectic vector
space $V$ and for every finite subgroup $G \subset \Aut V$ which
preserves the symplectic form, every crepant resolution of the
quotient singularity $X = V/G$ carries a canonical holomorphic
symplectic form $\Omega$. The form $\Omega$, which will be very
important for all our constructions, is described in this section.

Let $V$ be a complex vector space equipped with a non-degenerate
symplectic form $\Omega_V$. Assume given a finite group $G \subset
\Aut(V)$ of automorphisms of the vector space $V$ which preserves
the form $\Omega_V$. Let $X = V/G$ be the quotient variety. The
quotient map $\sigma:V \to X = V/G$ is \'etale over a smooth open
dense subset $X_0 \subset X$, moreover, it is a Galois covering with
Galois group $G$. Since the form $\Omega$ is $G$-invariant, it
defines a non-degenerate holomorphic $2$-form $\Omega \in
H^0(X_0,\Omega^2(X_0))$ over the subset $X_0 \subset X$. Moreover,
for every smooth resolution $\pi:\wt{X} \to X$ the pull-back
$\pi^*\Omega$ defines a holomorphic $2$-form over the preimage
$\wt{X}_0 = \pi^{-1}(X_0) \subset \wt{X}$.

\begin{defn}
A smooth resolution $\pi:\wt{X} \to X$ is called {\em symplectic} if
the canonical $2$-form $\pi^*\Omega$ on the open subset $\wt{X}_0
\subset \wt{X}$ extends to a non-degenerate symplectic form on the
whole smooth variety $\wt{X}$.
\end{defn}

Note that since we only consider irreducible resolutions $\wt{X}$,
the open subset $\wt{X}_0 \subset \wt{X}$ is always dense, and such
an extension $\Omega$ is unique.

The $n$-th power $\Omega^n$ of a non-degenerate symplectic $2$-form
on an $n$\--di\-men\-si\-onal smooth algebraic variety $Y$ gives a
section of the canonical bundle $K_Y$ without zeroes. Therefore
every symplectic resolution $\pi:\wt{X} \to X$ is crepant.

Our first result claims that the converse is also true. Namely, let
$\pi:\wt{X} \to X$ be a smooth crepant resolution. Then we have the
following.

\begin{prop}\label{crep=>sympl}
The pull-back $\pi^*\Omega$ of the holomorphic symplectic form
$\Omega$ on $X_0 \subset X$ extends to a non-degenerate symplectic
form on the whole smooth variety $\wt{X}$.
\end{prop}

\proof\footnote{This proof, as in fact everything in the paper,
works over an arbitrary algebraically closed field of characteristic
$0$. For a simpler proof which uses holomorphic geometry, see
\cite{Vnew} or \cite{Beau}.} The form $\pi^*\Omega$ is defined on
the open subset $\wt{X}_0 = \pi^{-1}(X_0) \subset \wt{X}$. A
$2$-form on an $n$-dimensional smooth algebraic variety is
non-degenerate if and only if the $n$-th power $\Omega^n$ has no
zeroes. Since the resolution $\pi:\wt{X} \to X$ is crepant, the
$n$-th power $\pi^*\Omega^n$ extends to the whole of $\wt{X}$ and
has no zeroes. Therefore it suffices to prove that the form
$\pi^*\Omega$ extends to a form $\Omega$ on the whole of $\wt{X}$.
Every such extension will be automatically non-degenerate.

We will prove slightly more.

\begin{lemma}\label{form.extends}
Let $V$ be an arbitrary connected smooth affine complex algebraic
variety equipped with an action of a finite group $G \subset
\Aut(V)$. Let $\pi:\wt{X} \to X$ be an arbitrary smooth resolution
of the quotient variety $X = V/G$. Denote by $\sigma:V \to X$ the
quotient map. 

Assume given a $k$-form $\alpha \in H^0(X_0,\Omega^k(X))$ on a
smooth dense open subset $X_0 \subset X$, and assume that the
pullback $\sigma^*\alpha$ extends to a form $\alpha_V \in
H^0(V,\Omega^k(V))$ on the whole smooth algebraic variety $V$.

Then the pullback $\pi^*\alpha$ extends to a form $\alpha \in
H^0(\wt{X},\Omega^k(\wt{X}))$ on the whole smooth algebraic variety
$\wt{X}$. 
\end{lemma}

\proof{} Since the variety $\wt{X}$ is smooth, it suffices to show
that $\pi^*\alpha$ extends to the complement to a closed subvariety
$Z \subset \wt{X}$ of codimension $\geq 2$. In other words, we have
to show that for every subvariety $Z \subset \wt{X} \setminus
\wt{X}_0$ of codimension $\codim Z = 1$, the form $\alpha$ extends
to an open neighborhood of a generic point $x \in Z$.

To do this, assume given an arbitrary (non-closed) point $x \in
\wt{X} \setminus \wt{X}_0$ with the residue field $K_x$ of
transcendental dimension $\dim \wt{X} - 1$ over $\C$. We have to
prove that the form $\pi^*\alpha$ extends to an open neighborhood of
the point $x$.

Let $\calo_{\wt{X},x}$ be the local ring of functions at $x$.  Since
$\wt{X}$ is smooth, the ring $\calo_{\wt{X},x}$ is a regular
discrete valuation ring and an algebra over the field $K_x$.
Consider the subset
$$
U = \Spec\calo_{\wt{X},x} \subset \wt{X}.
$$
The scheme $U$ is a smooth open curve over the field $K$. It has two
points: the special point $x \in U$, and the generic point $\eta \in
U$, which coincides with the generic point of the irreducible
variety $\wt{X}$. It suffices to prove that the form $\pi^*\alpha$
extends from the generic point $\eta \in U$ to the whole scheme $U$.

Since the map $\pi:\wt{X} \to X$ is generically one-to-one, it
induces an isomorphism between the generic point $\eta \in U \subset
\wt{X}$ and the generic point $\pi(\eta) \in X$. The quotient map
$\sigma:V \to X$ is a Galois covering over $\pi(\eta) \subset X$
with Galois group $G$.

Let $\wt{U}$ be the normalization of the curve $U$ in the field of
functions $K(V)$. This is a curve over the residue field $K_x$.

The curve $\wt{U}$ is normal, therefore it is smooth. Since the
curve $U$ is also smooth, the canonical projection $\tau:\wt{U} \to
U$ is flat. The group $G$ acts naturally on $\wt{U}$, and the map
$\tau:\wt{U} \to U$ induces an isomorphism $U \cong \wt{U}/G$. We
have a commutative square
$$
\begin{CD}
  \wt{U}     @>{\tau}>>        U     \\
@V{\pi}VV                  @VV{\pi}V \\
    V        @>{\sigma}>>      X
\end{CD}. 
$$
The pullback $\tau^*\pi^*\alpha$ is by definition a $G$-invariant
$k$-form on the generic subset $\tau^{-1}(\eta) \subset \wt{U}$. But
this form coincides with the pullback $\pi^*\sigma^*\alpha$, and the
pullback $\sigma^*\alpha$ by assumption extends to form $\alpha_V$
on the whole $V$. We conclude that the pullback $\tau^*\pi^*\alpha$
extends to a $G$-invariant $k$-form $\alpha$ on the whole curve
$\wt{U}$.

To finish the proof, it remains to apply the following simple
algebraic fact.

\begin{claim}
Let $k$ be a field of finite degree of transcendence over $\C$, let
$K_1$ be a Galois extension of the local field $K = k((x))$ with
Galois group $G$, and let $O_1$ be the integral closure of the
integer ring $O=k[[x]]$ in $K_1$. For every integer $n$, the have a
canonical isomorphism of $O$-modules
$$
\Omega^n(O/\C) \cong \Omega^n(O_1/\C)^G.
$$
(Here $\Omega^n$ is understood as the $n$-th exterior power of the
(flat) module $\Omega^1$ of K\"ahler differentials.)
\end{claim}

\proof{} If the extension $K_1/K$ is unramified, the claim is
obvious. Assume that the field $K_1$ is totally ramified over $K$,
so that the residue field $k_1$ of the local ring $O_1$ coincides
with $k$. Let $f \in O_1$ be the uniformizing element in the
discrete vaulation ring $O_1$. Since the field $k$ is of
characteristic $\cchar k = 0$, the Galois group $G$ is cyclic, say
$G = \Z/l\Z$. Moreover, $O_1 = k((f))$, $f^l \in O \subset O_1$
uniformizes $O$, and the module $\Omega^1(O_1/\C)$ of K\"ahler
differentials splits canonically into the direct sum
$$
\Omega^1(O_1/\C) = \left(\Omega^1(k/\C) \otimes_k O_1\right) \oplus
O_1 \cdot df. 
$$
Therefore the module $\Omega^n(O_1/\C)$ splits into the direct sum
$$
\Omega^n(O_1/\C) = \left(\vphantom{\Omega^{n-1}(k/\C)}
\Omega^n(k/\C) \otimes_k O_1 \right)
\oplus \left(\Omega^{n-1}(k/\C) \otimes_k O_1 \cdot df\right).
$$
The inclusion $\Omega^n(O/\C) \subset \Omega^n(O_1/\C)^G$ is
obvious. Moreover, since $O_1^G = O$ and $\Omega^i(k/\C)$ is
$G$-invariant for every $i$, we have the inverse inclusion
$$
\left(\Omega^n(k/\C) \otimes_k O_1 \right)^G = \Omega^n(k/\C)
\otimes_k O \subset \Omega^n(O/\C).
$$ 
It remains to prove the inclusion
$$
\left(\Omega^{n-1}(k/\C) \otimes_k O_1 \cdot df\right)^G 
\subset \Omega^{n-1}(O/\C).
$$
Since $O_1 = k[[f]]$, we have
$$
\Omega^{n-1}(k/\C) \otimes_k O_1 \cdot df = \prod_{p \geq 0}
\Omega^{n-1}(k/\C) \cdot f^pdf.
$$
To finish the proof, we note that a $1$-form $f^pdf$ is
$G$-invariant if and only if $p+1 = l(q+1)$ for some integer $q$,
and in this case we have
$$
\Omega^{n-1}(k/\C) \cdot f^pdf = \Omega^{n-1}(k/\C) \cdot x^qdx
\subset \Omega^{n-1}(O/\C). \qquad\qquad\qquad\square 
$$

Assume now given a smooth crepant resolution $\pi:\wt{X} \to X$ of a
symplectic quotient singularity $X = V/G$ which is equipped with the
holomorphic symplectic form $\Omega$ provided by
Proposition~\ref{crep=>sympl}. We will need the following
fundamental fact.

\begin{prop}\label{pullbacks}
Let $Y$ be a smooth algebraic variety equipped with maps $\pi_1:Y
\to V$ and $\pi_2:Y \to \wt{X}$ so that the square
\begin{equation}\label{pb}
\begin{CD}
Y            @>{\pi_2}>>   \wt{X}\\
@V{\pi_1}VV               @VV{\pi}V\\
V            @>{\sigma}>>    X
\end{CD}
\end{equation}
is commutative. Then we have
$$
\pi_1^*\Omega_V = \pi_2^*\Omega,
$$
where $\Omega_V$ is the canonical symplectic form on $V$, and
$\Omega$ is the canonical symplectic form on $\wt{X}$.
\end{prop}

\proof\footnote{A different proof based on the same general idea can
be found in \cite[Proposition~4.5]{Vhilb}.} Since $Y$ is smooth, it
suffices to prove the claim after replacing $Y$ with its generic
point $y \in Y$. Moreover, we can assume that $\pi_2(y) \not\in
\pi^{-1}(X_0) \subset \wt{X}$, since otherwise the claim follows
from definition.

Heuristically, the idea of the proof is to embed $y$ as a special
point $y \in Y$ into a smooth curve $Y \subset \wt{X}$ which
satisfies the conditions of the proposition, and whose generic point
lies in $\pi^{-1}(X_0) \subset \wt{X}$.

Let $x = \pi_2(y) \subset \wt{X}$ be the image of the point $y$. 
Consider the local ring $\calo_{\wt{X},x}$ of functions near $x$,
and let
$$
U = \Spec \calo_{\wt{X},x} \subset \wt{X}
$$ 
be the open neighborhood of $x \subset \wt{X}$. The intersection $U
\cap \wt{X}_0$ is not empty, since it contains the generic point of
$\wt{X}$. Take an arbitrary non-trivial function $f$ on $U$ which
vanishes on $U \cap (\wt{X} \setminus \wt{X}_0)$, and choose a point
$z \in U$ in the generic fiber of the map $f:U \to \A^1$ which is
closed in this generic fiber. Then the field of definition of the
point $z$ is of degree of transcedence $(\dim Y + 1)$ over $\C$, and
the closure $Z = \overline{z} \subset U$ is an open curve in $U$
with two points: the generic point $z \in Z$ and the special point
$y \in Z$.

Consider the pullback $\wt{Z} = Z \times_X V$, and let $\wh{Z}$ be
the normalization of the curve $\wt{Z}$. The curve $\wh{Z}$ is a
smooth open curve equipped with canonical maps $\pi_{Z,1}:\wh{Z} \to
V$ and $\pi_{Z,2}:\wh{Z} \to \wt{X}$ which satisfy the conditions of
the propositon. Moreover, the projection $\pi_2$ maps every generic
point of the curve $\wh{Z}$ into the open subset $\wt{X}_0 \subset
\wt{X}$. Therefore we know that
\begin{equation}\label{Z.pb}
\pi_{Z,1}^*\Omega_V = \pi_{Z,2}^*\Omega.
\end{equation}
Now, since the square \eqref{pb} is commutative, the maps $\pi_1:y
\to V$, $\pi_2:y \to Z \subset \wt{X}$ both factor through a map
$\nu:y \to \wt{Z} = V \times_X Z$. Therefore we can replace $y$ with
its image $\wt{y} = \nu(y) \subset \wt{Z}$, which is a closed point
in $\wt{Z}$.

Choose a closed point $\wh{y} \in \wh{Z}$ \'etale over $\wt{y}
\subset \wt{Z}$, and let $\tau:\wh{y} \to \wt{y}$ be the
projection. By construction we have $\pi_{Z,1}|_{\wh{y}} = \pi_1
\circ \tau$ and $\pi_{Z,2}|_{\wh{y}} = \pi_2 \circ \tau$, so that
\eqref{Z.pb} yields
$$
\tau^*\pi_1^*\Omega_V = \tau^*\pi_2^*\Omega.
$$
Since $\tau:\wh{y} \to \wt{y}$ is \'etale, this implies the claim.
\endproof

We end this section with the following corollary of
Lemma~\ref{form.extends}, which is a reformulation of
Theorem~\ref{equivar}~\thetag{i}.

\begin{corr}\label{field.lifts}
Let $\TT_V$ be an arbitrary $G$-invariant vector field on the
symplectic complex vector space $V$, and let $\TT_0$ be the induced
vector field on the open subset $X_0 \subset X$ of the quotient
variety $X = V/G$.

Then for every smooth crepant resolution $\pi:\wt{X} \to X$, the
vector field $\TT_0$ lifts to a vector field $\TT$ on the whole
smooth algebraic variety $X$.
\end{corr}

\proof{} The holomorphic symplectic form $\Omega$ induces an
isomorphism 
$$
\Omega^1(X_0) \cong \T(X_0)
$$ 
between the sheaf $\Omega^1(X_0)$ of $1$-forms and the sheaf
$\T(X_0)$ of vector fields, and a compatible isomorphism
$\Omega^1(\wt{X}) \cong \T(\wt{X})$. Let $\alpha \in
H^0(X_0,\Omega^1(X_0))$ be the $1$-form corresponding to the vector
field $\TT_0 \in H^0(X_0,\T(X_0))$. It suffices to prove that the
pullback $\pi^*\alpha$ extends to the whole $\wt{X}$. This follows
from the assumption and from Lemma~\ref{form.extends}.  \endproof

\begin{remark}\label{CalabiYau}
In fact, an analogous result holds in the more general situation of
Calabi-Yau manifolds. More precisely, instead of symplectic vector
spaces one can consider vector spaces equipped with a volume
form. If the group $G$ preserves the volume form, then
Corollary~\ref{field.lifts} immediately extends to crepant
resolutions of the quotient $X = V/G$. The proof is the same, with
one modification -- one has to identify vector fields with
$(n-1)$-forms rather than with $1$-forms.
\end{remark}

\section{Stratification by rank}

As in the last section, let $V$ be a complex vector space equipped
with a non-degenerate symplectic form $\Omega_V$, and let $G \subset
\Aut(V)$ be a finite group of automorphisms of the vector space $V$
which preserves the form $\Omega_V$. We will now introduce a
canonical stratification on the singular quotient variety $X = V/G$,
which we will call {\em stratification by rank}.

To do this, for every vector $v \in V$ let $G_v \subset G$ be the
stabilizer subgroup of the vector $v$, that is, the subgroup of
elements $g \in G$ such that $g \cdot v = v$. Define the {\em rank}
$\rk v$ of the vector $v \in V$ as one half of the codimension
$$
\rk v = \frac{1}{2}\left(\dim V - \dim V^{G_v}\right)
$$
of the subspace $V^{G_v} \subset V$ of $G_v$-invariant vectors. The
number $\rk v$ obviously depends only on the image $\sigma(v) \in X$
of the vector $v \in V$ under the quotient map $\sigma:V \to X =
V/G$.

Note that the restriction of the symplectic form $\Omega_v$ to the
subspace $V \subset V^{G_v}$ is necessarily
non-degenerate. Therefore $\dim V^{G_v}$ is always even, and $\rk v$
is an integer for every vector $v \in V$.

For every integer $k \geq 0$, let $X_k \subset X$ be the subset of
points $x \in X$ such that the corresponding $G$-orbit in $V$
consists of vectors of rank $k$. The subsets $X_k$ are locally
closed subvarieties in $X$, and $X_p \subset \overline{X_q}$ implies
$p \geq q$. The largest subset $X_0 \subset X$ is open and
dense. The same is true for the unions
$$
X_{\leq p} = \bigcup_{q \leq p}X_q.
$$

\begin{lemma}\label{strat}
The subset $X_k \subset X$ is a (non-connected) smooth algebraic
variety of dimension $\dim Y_k = \dim X - 2k$. Moreover, the
projection map $\sigma^{-1}(X_k) \to X_k$ is \'etale.
\end{lemma}

\proof{} Let $V_k = \sigma^{-1}(X_k) \subset V$. Consider a vector $v
\in V_k$ and let $x = \sigma(v) \in X_k$. Let $\wh{V}$ be the formal
completion of $V$ at $v$, and let $\wh{X}$ be the formal completion
of $X$ at $x$. The projection map $\sigma:\wh{V} \to \wh{X}$ induces
an isomorphism $\wh{V}/G_v \cong \wh{X}$. Moreover, it identifies
the formal completion $\wh{X}_k \subset \wh{X}$ of the subvariety
$X_k \subset X$ at $x$ with the quotient $\wh{V}_k/G_v$ of formal
completion $\wh{V}_k \subset \wh{V}$ of the subvariety $V_k \subset
V$ at $v \in V_k$.

But the completion $\wh{V}_k$ is isomorphic to the completion at $0$
of the vector subspace $V^{G_v} \subset V$ of $G_v$-invariant
vectors. Therefore $G_v$ acts trivially on $\wh{V}_k$, the quotient
$\wh{V}_k/G_v = \wh{V}_k$ is smooth and the quotient map $\sigma:V_k
\to X_k$ is \'etale at $v$. 

The dimension formula follows directly from the definition of the
rank $\rk v$.
\endproof

As a consequence of this lemma, we see that the $G$-invariant
symplectic form $\Omega$ on the vector space $V$ defines a $2$-form
$\Omega_k$ on every stratum $X_k$. Since the restriction
$\Omega|_{V^{G_v}}$ is non-degenerate for every $G_v \subset G$, all
these forms are non-degenerate holomorphic symplectic forms.

Let $Y \subset X_k$ be a connected component of the locally closed
subvariety $X_k \subset X$. While for $k \geq 1$ the subvariety $Y
\subset X$ lies in the singular locus of $X$, it still admits a sort
of a ``tubular neighborhood''. Namely, let $\wh{Y}$ be the formal
completion of $X$ along $Y \subset X$. Moreover, choose a connected
component $V_Y \subset \sigma^{-1}(Y)$ of the preimage
$\sigma^{-1}(Y) \subset V$ of $Y \subset X$ under the quotient map
$\sigma:V \to X$, and let $\wh{V}_Y$ be the formal completion of $V$
along $V_Y \subset V$.

\begin{lemma}\label{tubular}
\begin{enumerate}
\item The stabilizer subgroup $G_v \subset G$ is the same for every
vector $v \in V_Y$. 

\item The quotient map $\sigma:V \to X$ induces an \'etale map
$$
\wh{V}_Y/G_v \to \wh{Y}.
$$

\item Let $V' \subset V$ be the $G_v$-invariant complement to the
subspace $V^{G_v} \subset V$ of vectors fixed by $G_v$, and let $X'
= V'/G_v$ be the quotient variety. There exists a canonical direct
product decomposition
$$
\wh{V}_Y \cong \wh{X'} \times V_Y, 
$$
where $\wh{X'}$ is the formal completion of $X'$ along $0 \in X'$.
\end{enumerate}
\end{lemma}

\proof{} \thetag{i} and \thetag{ii} are clear. The direct product
decomposition in \thetag{iii} is induced by the the direct sum
decomposition $V = V' \oplus V^{G_v}$. \endproof

\comment
We will be especially interested in connected components of the
subset $X_1 \subset X$. We call these components {\em mirrors}, for
the following reason.

\begin{lemma}\label{mirrors}
For every vector $v \in V$ of rank $1$ the stabilizer group $G_v
\subset G$ is generated by a single reflection $g \in G$ of order
$2$. The images 
$$
\sigma(v_1),\sigma(v_2) \in X_1
$$ 
of any two vectors $v_1, v_2 \in V$ of rank one lie in the same
connected component $Y \subset X_1$ if and only if the corresponding
reflections $g_1,g_2 \in G$ are conjugate within $G$.
\end{lemma}

\proof{} Since the subgroup $G_v$ acts trivially on a subspace in $V$
of codimension $1$, it must lie in the group $\R^*$ of invertible
real numbers, and since $G_v$ is finite, it must coincide with
$\{\pm 1\} \subset \R^*$. The second claim is clear.
\endproof
\endcomment

It is easy to give a combinatorial description of the set of
components $Y$ of every stratum $X$, but we will do it only for the
simplest case and under the restrictive assumptions of
Theorem~\ref{unique}. It is the only case which we will need later
on.

\begin{lemma}\label{mirrors}
Let $V_\R$ be a real vector vector space equipped with an action of
a finite group $G$. Let $V_o = V_\R \otimes \C$, let $V = V_o \oplus
V_o$, and let $X = V/G$ be the quotient of the vector space $V$ by
the natural $G$-action.

Then the connected components $Y \subset X$ os the first stratum
$X_1 \subset X$ of the stratification by rank are in a natural
one-to-one correspondence with the conjugacy classes of reflections
in the subgroup $G \subset V_o$. Moreover, every reflection in $G
\subset V_o$ is of order $2$.
\end{lemma}

\proof{} First of all, every reflection $g \in G \subset \Aut V_o$
preserves the real subspace $V_\R$. Therefore it acts by
multiplication by a {\em real} scalar, that is, $\pm 1$. Thus $g$ is
of order two.

Let $g \in G \subset \Aut V_o$ be a reflection. Consider the
subspace $V^g \subset V$ of $g$-invariant vectors and the image
$\sigma(V^g) \subset X$ under the quotient map $\sigma:V \to X$. The
generic vector $v \in V^g$ does not lie in any smaller $G$-invariant
subspace in $V$. Therefore $\rk v = 1$, and the intersection $Y^g =
\sigma(V^g) \cap X_1$ is dense in $\sigma(V^g)$. 

Since $\sigma(V^g)$ is irreducible, the intersection $Y^g$ is a
connected component of the stratum $X_1$. Moreover, for two
reflections $g_1,g_2 \in G$ we have $Y^{g_1} = Y^{g_2} \subset X_1$
if and only if the reflections $g_1$ and $g_2$ are conjugate within
$G$. This shows that the correspondence $g \mapsto Y^g$ is
injective. 

Finally, for every connected component $Y \subset X_1$ let $v \in V$
be an arbitrary vector in the preimage $\sigma^{-1}(Y) \subset
V$. Then $\rk v = \codim V_\R^{G_v} = 1$, and the stabilizer
subgroup $G_v \subset G$ is in fact a subgroup in the group
$\Aut(V_\R/V_\R^{G_v}) = \R^*$. Therefore $G_v$ contains a single
non-trivial elementr $g \neq 1 \in G_v = G$, the element $g \in G
\subset \Aut V_o$ is a reflection, and we have $Y = Y^g$. Since $Y$
is arbitrary, this implies that the correspondence $g \mapsto Y^g$
is one-to-one.  \endproof

Assume now given a smooth crepant resolution $\pi:\wt{X} \to X$ of
the quotient variety $X$. The stratification $X_k$ induces a locally
closed stratification
$$
\wt{X}_k = \pi^{-1}(X_k) \subset \wt{X}.
$$
The strata $\wt{X}_k$ are no longer necessarily smooth. 

As in \cite[Proposition 4.16]{Vhilb}, Proposition~\ref{pullbacks}
immediately implies the following.

\begin{prop}\label{semismall}
For every stratum $\wt{X}_k \subset \wt{X}$ we have
$$
\dim\wt{X}_k \leq \dim_X - k.
$$
\end{prop}

\proof{} Let $Y \subset \wt{X}_k$ be the non-singular part of the
stratum $X_k$. Let $y \in Y$ be an arbitrary point, let $T_yY
\subset T_y\wt{X}$ be the tangent space to $Y$ at the point $y$, and
let $T_{vert} \subset T_yY$ be the kernel of the differential
$d\pi:T_yY \to T_{\pi(y)}X_k$. By Proposition~\ref{pullbacks}, the
restriction of the canonical non-degenerate $2$-form $\Omega$ to the
subvariety $Y \subset \wt{X}$ coincides with the pullback of the
canonical non-degenerate $2$-form on the stratum $X_k$. Therefore
$T_{vert}$ is orthogonal to $T_yY$ with respect to the form
$\Omega$. This implies that
$$
\dim T_{vert} + \dim T_yY \leq \dim T_y\wt{X}.
$$
But $\dim T_yY = \dim Y$, $\dim T_y\wt{X} = \dim X$ and $\dim
T_{vert} = \dim Y - \dim X_k$. Therefore
$$
2 \dim Y \leq \dim X + \dim X_k = 2 \dim X - 2k,
$$
which proves the claim.
\endproof

In particular, we see that the union 
$$
\wt{X}_{\leq 1} = \wt{X}_0 \cup \wt{X}_1 \subset \wt{X}
$$
is an open subset whose complement is of codimension $> 1$. This
fact will be crucial for all our constructions.

\begin{remark}
Later (Proposition~\ref{tubular.res}) we shall see that for every
$k$ the stratum $\wt{X}_k \subset \wt{X}$ is equidimensional over
the stratum $X_k \subset X$. Thus Proposition~\ref{semismall} means
that the resolution $\pi:\wt{X} \to X$ is semismall.
\end{remark}

\section{The Picard group}\label{sec.pic}

We will now use the rank stratification to study resolutions of the
symplectic quotient variety $X = V/G$. We begin with the following
general observation.

\begin{lemma}\label{blowup}
The Picard group $\Pic(X_0)$ of the non-singular open stratum $X_0
\subset X$ is a torsion group. Therefore, every projective
resolution $\pi:\wt{X} \to X$ is canonically isomorphic to the
blow-up
$$
\Bl\left(X,\E\right) \to X
$$
of a sheaf $\E \subset \calo_X$ of ideals on $X$.
\end{lemma}

\proof{} The quotient map $\sigma:V \to X$ is \'etale over $X_0
\subset X$, so that modulo torsion $\Pic(X_0)$ is a subgroup in $V_0
= \sigma^{-1}(X_0) \subset V$. But the complement to $V_0$ in $V$ is
codimension $> 1$ by Lemma~\ref{strat}. Therefore $\Pic(V_0) =
\Pic(V) = 0$. This proves the first claim.

To prove the second claim, let $L$ be a line bundle on $\wt{X}$
which is very ample for the map $\pi$, so that
$$
\wt{X} = \Proj\bigoplus_{k \geq 0}\pi_*(L^k).
$$
The projective resolution $\pi:\wt{X} \to X$ is by definition
one-to-one over an open dense subset $U \subset X$. Since $X$ is
normal, we can assume that the complement $X \setminus U \subset X$
is of codimension $\codim X \setminus U \geq 2$. 

Replace $U$ with $U \cap X_0$, so that $\Pic(U) = \Pic(X_0)$, and
denote by $\wt{U} = \pi^{-1}(U) \subset \wt{X}$ the preimage of the
open subset $U \subset X$. Note that since the subset $X \setminus
X_0 \subset X$ is of codimension $\codim X \setminus X_0 \geq 2$, we
still have $\codim X \setminus U \geq 2$.

By assumption $\wt{U} \cong U$ and $\Pic(\wt{U}) = \Pic(U) =
\Pic(X_0)$ is a torsion group. Replacing the line bundle $L$ with
its positive power, we can assume that it is trivial on $\wt{U}$.

Choose a trivialization map $L \to \calo_{\wt{U}}$, or,
equivalently, $\pi_*L \to \calo_{U}$. Since the complement to $U$ in
the normal variety $X$ is of codimension $\geq 2$, this
trivialization map extends to an embedding
$$
\pi_*L \to \calo_X.
$$
Denoting $\E = \pi_*L \subset \calo_X$, we get the second claim.
\endproof

Let now $\pi:\wt{X} \to X$ be a smooth projective crepant resolution
of the quotient variety $X = V/G$. Lemma~\ref{blowup} allows us to
prove the first of the results announced in Section~1.

\proof[Proof of Theorem~\ref{equivar}.] \thetag{i} is
Corollary~\ref{field.lifts}. To prove \thetag{ii}, note that the Lie
algebra $\hh$ of the group $H$ acts naturally on the resolution
$\wt{X}$ by \thetag{i}. Therefore the ideal $\E \subset \calo_X$
provided by Lemma~\ref{blowup} is $\hh$-invariant. Since the group
$H$ is by assumption connected, the ideal $\E \subset \calo_x$ is
also $H$-invariant, which implies that the $H$-acton lifts to the
resolution $\wt{X} = \Bl(X,\E)$.

To derive \thetag{iii} from \thetag{ii}, it suffices to notice that
the algebraic group $GL(V,G)$ of $G$-equivariant linear
automorphisms of the vector space $V$ is connected.  \endproof

Next, we note that all the results used to establish
Theorem~\ref{equivar}, in particular, Lemma~\ref{blowup} and
Corollary~\ref{field.lifts}, carry over literally to the formal
setting of Theorem~\ref{formal}, with the same proofs. We will use
this to show that Theorem~\ref{formal} is also a corollary of
Lemma~\ref{blowup}.

\proof[Proof of Theorem~{\normalfont\ref{formal}}.]  Assume given a
projective smooth formal crepant resolution $\wh{\pi}:\X \to \wh{X}$
of the completion of the quotient variety $X = V/G$ at $0 \in
X$. Let $\wh{\E} \subset \wh{\calo}_X$ be the ideal such that $\X =
\Bl(\wh{X},\wh{\E})$. 

Let the group $\C^*$ act on the vector field $V$ by homoteties. This
action commutes with the $G$-action and defines a $\C^*$-action on
the quotient $V/G$ or, equivalently, a grading on the algebra
$\calo_X$ of functions on $X$. This grading induces a decreasing
filtration on the completion $\wh{\calo}_X$, and the associated
graded quotient with respect to this filtration is the algebra
$\calo_X$. Let $\E \subset \calo_X$ be the associated graded
quotient to the ideal $\wh{\E} \subset \wh{\calo}_X$.

The differential of the $\C^*$-action on $V$, namely, the Euler
vector field $\xi_V$ on $V$, induces a vector field $\xi$ on the
quotient $X = V/G$ and on the completion $\wh{X}$.  By
Corollary~\ref{field.lifts}, the vector field $\xi$ preserves the
ideal $\wh{\E} \subset \wh{\calo}_X$. Therefore this ideal is
isomorphic to the completion of its associated graded quotient $\E
\subset \calo_X$. Let
$$
\wt{X} = \Bl\left(X,\E\right)
$$
be the blow-up of the ideal $\E \subset \calo_X$, and let
$\pi:\wt{X} \to X$ be the projection. 

By construction the map $\pi$ is projective. Moreover, since
$\wh{\E} \subset \wh{\calo}_X$ is the completion of $\E \subset
\calo_X$, the completion of $\pi:\wt{X} \to X$ along $\pi^{-1}(0)
\subset X$ coincides with $\wh{\pi}:\X \to X$. Therefore $\wt{X}$ is
smooth in an open neighborhood $U \subset X$ of $\pi^{-1}(0)$, and
the canonical bundle $K_U$ is trivial. 

But the ideal $\E \subset \calo_X$ is by construction homogenous
with respect to the grading given by the Euler vector
field. Therefore the $\C^*$-action on $X = V/G$ lifts to a
$\C^*$-action on $\wt{X}$.  Since for every point $x \in \wt{X}$ we
have
$$
\lim_{\begin{Sb}\lambda \to 0\\\lambda \in \C^*\end{Sb}} \pi(\lambda
\cdot x) = 0,
$$
we have $\lambda \cdot x \in U \subset \wt{X}$ for some $\lambda \in
\C^*$, and this implies that $\wt{X}$ is smooth
everywhere. Moreover, since the canonical bundle $K_{\wt{X}}$ is
$\C^*$-equivariant, and the canonical bundle $K_U$ is trivial, the
whole bundle $K_{\wt{X}}$ is trivial, and the resolution $\pi:\wt{X}
\to X$ is crepant.  \endproof

Consider now the rank stratification of the variety $X$, and let $Y
\subset X_k$ be a connected component of an arbitrary stratum $X_k
\subset X$. We will use Lemma~\ref{blowup} to describe the structure
of the resolution $\pi:\wt{X} \to X$ near the subvariety $Y \subset
X$.

To do this, return to the setting of Lemma~\ref{tubular}. Denote by
$\wh{Y}$ the formal completion of $X$ near the subvariety $Y \subset
X$. Choose a connected component $V_Y \subset \sigma^{-1}(Y)$ of the
preimage $\sigma^{-1}(Y) \subset V$ of the subvariety $Y \subset X$
under the quotient map $V \to X$, and let $\wh{V}_Y$ be the
completion of $V$ along $V_Y \subset V$.

Choose an arbitrary vector $v \in V_Y$, let $G_v \subset G$ be the
stabilizer subgroup of the vector $v$, and let $V' \subset V$ be the
unique $G_v$-invariant complement to the subspace $V^{G_v} \subset
V$ of $G_v$-invariant vectors.

By Lemma~\ref{tubular}, we have a canonical \'etale map
$$
\sigma:\wh{V}_Y\to \wh{Y}.
$$
Moreover, the quotient $X' = V'/G_v$ does not depend on the
choice of the point $v \in V_Y$, and we have a direct product
decomposition 
$$
\wh{V}_Y \cong \wh{X'} \times V_Y,
$$
where $\wh{X'}$ is the completion of the quotient $X'$ near $0 \in
X'$.

Let $\Y_0$ be the completion of $\wt{X}$ along the subvariety
$\wt{Y} = \pi^{-1}(Y) \subset \wt{X}$, and let $\Y$ be the fibered
product given by the diagram
$$
\begin{CD}
\Y              @>>>         \Y_0\\
@V{\pi}VV                  @VV{\pi}V\\
\wh{V}_Y/G_v @>{\sigma}>>    \wh{Y}.
\end{CD}
$$

\begin{prop}\label{tubular.res}
There exists a a formal smooth projective crepant resolution
$\pi':\X' \to \wh{X'}$ of the quotient variety $X' = V'/G_v$ and a
direct product decomposition
\begin{equation}\label{prod}
\Y \cong \X' \times V_Y,
\end{equation}
such that 
$$
\pi = \pi' \times \id:\Y \cong \X' \times V_Y \to \wh{V}_Y/G_v \cong
\wh{X'} \times V_Y. 
$$
\end{prop}

\proof{} Since the map $\sigma:V_Y \times \wh{X'} \to \wh{Y}$ is
\'etale, Lemma~\ref{blowup} implies that
$$
\Y = \Bl\left(V_Y \times \wh{X'},\E\right),
$$
where $\E$ is a sheaf of ideals on $V_Y \times \wh{X'}$. Let
$\pi_1:\wh{V}_Y \cong \wh{X'} \times V_Y \to \wh{X'}$ be the
projection onto the first factor. To obtain the direct product
decomposition \eqref{prod}, it suffices to prove that $\E =
\pi_1^*\E'$ for some sheaf of ideals $\E'$ on $\wh{X'}$.

Every vector $v \in V$ defines a constant vector field $\TT_v$ on
the vector space $V$. If the vector $v \in V^{G_v}$ is
$G_v$-invariant, then the vector field $\TT_v$ is also
$G_v$-invariant and descends to a vector field on the quotient
$V/G_v$.  Completing along $V_Y \subset V/G_v$, we obtain a vector
field $\TT_v$ on $\wh{V}_Y \cong \wh{X'} \times V_Y$.

Vector fields $\TT_v$, $v \in V^{G_v}$ are parallel to the fibers of
the projection $\pi_1:\wh{V}_Y/G_V \to \wh{X'}$ and generate the
relative tangent bundle of $\wh{V}_Y/G_v$ over $\wh{X'}$. Therefore,
to prove that the sheaf of ideals $\E$ satifies $\E = \pi_1^*\E'$
for some $\E'$ on $\wh{X'}$, it suffices to prove that $\E$ is
preserved by the vector field $\TT_v$ for an arbitrary $v \in
V^{G_v}$. This is equivalent to proving that the vector field
$\TT_v$ lifts to a vector filed on the resolution $\Y$ of the
variety $\wh{X'} \times V_Y$. But this follows from
Lemma~\ref{form.extends} by exactly the same argument as
Corollary~\ref{field.lifts}. Thus the ideal $\E$ on $\wh{X'} \times
V_Y$ is preserved by all the vector fields $\TT_v$, so that we have
$\E = \pi_1^*\E'$ for some ideal $\E'$ on $\wh{X'}$. 

It remains to prove that the resolution
$$
\pi':\X' = \Bl\left(\wh{X'},\E'\right) \to \wh{X'}
$$
is smooth and crepant. But this immediately follows from the
corresponding properties of the resolution $\pi:\Y \to \wh{X'}
\times V_Y$ and from the direct product decomposition \eqref{prod}.
\endproof

Our last general result, Theorem~\ref{inductive}, follows
immediately from Theorem~\ref{formal} and
Proposition~\ref{tubular.res}.

We now restrict the generality and introduce the same assumptions on
$\langle V,G \rangle$ as in Theorem~\ref{unique}. Namely, assume
that the symplectic vector space $V$ is of the form $V = V_o \oplus
V_o^*$, where $V_o \subset V$ is a $G$-invariant Lagrangian
subspace. Assume additionally that we have $V_o = V_\R \otimes_\R
\C$ for some real vector space $V_\R$, and that the $G$-action on
$V_o$ is induced by a $G$-action on $V_\R$. 

Under these assumptions, and in the case of strata of small rank,
Proposition~\ref{tubular.res} can be made more precise.

\begin{corr}\label{0.and.1}
Let $\pi:\wt{X} \to X$ be  smooth projective crepant resolution of
the quotient variety $X = (V_o \oplus V_o^*)/G$.
\begin{enumerate}
\item Over the open stratum $X_0 \subset X$, the map $\pi:\wt{X} \to
X$ is one-to-one.
\item Let $\wh{X_1}$ be the completion of $X$ along the stratum $X_1
\subset X$ of codimension $2$, and let $\X_1$ be completion of
$\wt{X}$ along the preimage $\wt{X}_1 = \pi^{-1}(X_1) \subset
\wt{X}$. Then the resolution $\pi:\X_1 \to \wh{X}_1$ is isomorphic
to the blow-up of the completion $\wh{X_1}$ along the closed
subvarety $X_1 \subset \wh{X_1}$.
\end{enumerate}
\end{corr}

\proof{} \thetag{i} is immediate from
Proposition~\ref{tubular.res}. \thetag{ii} follows from
Proposition~\ref{tubular.res} and Proposition~\ref{dim.1}.
\endproof

To combine these two particular cases, let
$$
X_{\leq 1} = X_0 \cup X_1 \subset X
$$
be the open dense subset of point of rank $\leq 1$ in $X$, and let
$\wt{X}_{\leq 1} \subset \wt{X}$ be its preimage in $\wt{X}$. 

\begin{corr}\label{leq.1}
The restriction $\pi:\wt{X}_{\leq 1} \to X_{\leq 1}$ of the
resolution $\pi:\wt{X} \to X$ to $\wt{X}_{\leq 1}$ is canonically
isomorphic to the blow-up
$$
\Bl\left(X_{\leq 1},X_1\right) \to X_{\leq 1}
$$
of the variety $X_{\leq 1}$ in the closed subvariety $X_1 \subset
X_{\leq 1}$.
\end{corr}

\proof{} Lemma~\ref{blowup} and Corollary~\ref{0.and.1}.
\endproof

Let now $M$ be the set of all connected components of the subvariety
$X_1 \subset X$, and consider the free $\Q$-vector space $\Q[M]$
generated by the set $M$. (Note that by Lemma~\ref{mirrors} the set
$M$ coincides with the set of conjugacy classes of reflections in
$M$.) Denote also by $\Q_+[M] \subset \Q[M]$ the subset of linear
combinations of elements of the set $M$ with positive coefficients.

\begin{lemma}
The rational Picard group $\Pic(\wt{X}) \otimes \Q$ of the smooth
resolution $\wt{X}$ is canonically isomorphic to $\Q[M]$, and the
isomorphism can be chosen in such a way that a class $[L] \in \Q[M]$
of a line bundle $L$ on $\wt{X}$ which is very ample for the map
$\pi:\wt{X} \to X$ lies in $\Q_+[M] \subset \Q[M]$.
\end{lemma}

\proof{} By Proposition~\ref{semismall}, the complement $\wt{X}
\setminus \wt{X}_{\leq 1}$ is of codimension $\geq 2$. Therefore we
have $\Pic(\wt{X}) \cong \Pic(\wt{X}_1)$. By
Corollary~\ref{0.and.1}~\thetag{ii} and Proposition~\ref{dim.1}, the
preimage $\wt{Y} = \pi^{-1}(Y) \subset \wt{X}_{\leq 1}$ of every
connected component $Y \subset X_1$ of the subvariety $X_1 \subset
X_{\leq 1}$ is a smooth divisor in $X_{\leq 1}$ isomorphic to
$$
\wt{Y} \cong Y \times \cp^1.
$$
Now, the preimage $\sigma^{-1}(Y) \subset V$ is an open subset in
the union 
$$
\bigcup_g V^g \subset V
$$
of fixed-points subspaces for reflections $g \subset G$ in the
conjugacy class corresponding to $Y$. We can replace this union with
a single space $V^g$ and obtain an isomorphism
$$
\overline{Y} \cong V^g/N(g) \subset X,
$$
where $N(g) \subset G$ is the normalizer subgroup of the subspace
$V^g \subset V$. The complement to $\sigma^{-1}(Y) \cap V^g$ in
$V^g$ is of codimension $\geq 2$, therefore $\Pic(\sigma^{-1}(Y)
\cap V^g) = 0$. Since the quotient map $\sigma:V \to X$ is \'etale
over $Y$ by Lemma~\ref{strat}, this implies that $\Pic(Y) \otimes \Q
= 0$. Therefore $\Pic(\wt{Y}) \otimes \Q = \Pic(\cp^1) \otimes \Q =
\Q$, and the restriction to $\wt{X}_1 \subset \wt{X}_{\leq 1}$
defines a map
$$
\res:\Pic(\wt{X}) \otimes \Q \cong \Pic(\wt{X}_{\leq 1}) \otimes \Q
\to \bigoplus_{Y \in M}\Pic(\wt{Y}) \otimes \Q = \Q[M].
$$
The map $\res$ is injective. Indeed, for every line bundle $L$ on
$\wt{X}_{\leq 1}$ with trivial restriction $L|_{\wt{X}_1}$ we have
$L \cong \pi^*\pi_*L$, which means that the sheaf $\pi_*L$ is a line
bundle. This line bundle is in turn isomorphic to
$$
\pi_*L \cong j_*j^*\pi_*L,
$$
where $j:X_0 \hookrightarrow X_{\leq 1}$ is the embedding map of the
open subset $X_0$. Since $\Pic(X_0)$ is torsion, we can assume that
$j^*\pi_*L \cong \calo_{X_0}$ is trivial. Therefore $\pi_*L =
j_*\calo_{X_0} = \calo_{\wt{X}_1}$ is trivial as well.

But the map $\res$ is also surjective. Indeed, the correspondence $Y
\mapsto \calo(\wt{Y})$ defines a map $\cl:\Q[M] \to
\Pic(\wt{X}_{\leq 1})$, and the composition
$$
\res \circ \cl:\Q[M] \to \Pic(\wt{X}_{\leq 1}) \to \Q[M]
$$
is multiplication by $2$.

To prove the last claim of the lemma, it suffices to choose for
every $Y \in M$ a relatively very ample line bundles as the
generator for the group $\Pic(\wt{Y}) \cong \Q$. To fix the
isomorphism once and for all, we will take for such a generator the
bundle $\calo(1)$ on $\wt{Y} \cong Y \times \cp^1$.  \endproof

We can now formulate and prove the main result of this section. For
every connected component $Y \in M$ of the subvariety $X_1 \subset
X$ denote by
$$
\m_Y \subset \calo_X
$$
the ideal of the closed subvariety $\overline{Y} \subset X$, and for
every element $l \in \Q_+[M]$ which is a linear combination of
generators with positive even integer coefficients,
\begin{equation}\label{even}
l = \sum_{Y \in M} 2k_Y [Y] \in \Q[M],
\end{equation}
denote by $\E_l \subset \calo_X$ the intersection ideal 
\begin{equation}\label{ideal}
\E_l = \bigcap_{Y \in M} \m_Y^{k_Y}.
\end{equation}

\begin{prop}\label{class}
Let $\pi:\wt{X} \to X$ be a smooth projective crepant resolution,
and let $L$ be a line bundle on $X$ which is very ample for the map
$\pi$. Assume that the class $\cl(L) \subset \Q[M]$ is of the form
\eqref{even}. Then the given resolution $\pi:\wt{X} \to X$ is
canonically isomorphic to the blow-up
$$
\Bl\left(X,\E_l\right) \to X
$$
of the ideal $\E_l \subset \calo_X$. In particular, if for two
smooth crepant resolutions $\pi:\wt{X} \to X$, $\pi':\wt{X}' \to X$
with relatively very ample line bundles $L$, $L'$ the classes
$\cl(L),\cl(L') \subset \Q[M]$ are proportional, then the given
resolutions $\pi:\wt{X} \to X$, $\pi':\wt{X}' \to X$ are isomorphic.
\end{prop}

\proof{} By Lemma~\ref{blowup} the resolution $\pi:\wt{X} \to X$ is
isomorphic to the blow-up of the ideal $\E = \pi_*L \subset
\calo_X$. Let $j_1:X_{\leq 1} \hookrightarrow X$ be the embedding
map of the open subset $X_{\leq 1} \subset X$, and let
$\wt{\jmath}_1:\wt{X}_{\leq 1} \hookrightarrow \wt{X}$ be the embedding
map of the open subset $\wt{X}_{\leq 1} \subset \wt{X}$. Since the
complement $\wt{X} \setminus \wt{X}_{\leq 1}$ is of codimension $>
1$, we have
$$
L \cong \wt{\jmath}_{1*}\wt{\jmath}_1^*L.
$$
Therefore
$$
\E = \pi_*L = \pi_*\wt{\jmath}_{1*}\wt{\jmath}_1^*L =
j_{1*}\pi_*\wt{\jmath}^*_1L.
$$
Since $\E_l \cong j_{1*}j^*_1\E_l$, it suffices to prove that $\E_l
= \pi_*L \subset \calo_{X_{\leq 1}}$ on $X_{\leq 1}$. This follows
immediately from Corollary~\ref{leq.1}.

To prove the second claim, it suffices to notice that both classes
$\cl(L)$ and $\cl(L')$ can be made equal and of the form
\eqref{even} by multiplication by an appropriate positive
integer. This corresponds to replacing $L$ and $L'$ with positive
powers, which does not change the resolutions $\pi:\wt{X} \to X$ and
$\pi':\wt{X}' \to X$.
\endproof

This immediately yields the following corollary, which is a
reformulation of Theorem~\ref{unique}.

\begin{corr}\label{unicor}
If the set $M$ consists of a single element, then every two
projective smooth crepant resolutions $\pi:\wt{X} \to X$,
$\pi':\wt{X}' \to X$ of the quotient variety $X$ are
isomorphic.\endproof
\end{corr}

\begin{remark}\label{semism.unique}
As we can see from the proof, it is not necessary to require
smoothness of the resolutions $\wt{X}$, $\wt{X}'$ in
Proposition~\ref{class} and Corollary~\ref{unicor}. It suffices to
require the following:
\begin{enumerate}
\item The varieties $\wt{X}$, $\wt{X}'$ are normal.
\item Over $X_{\leq 1} \subset X$, we have 
$$
\wt{X}_{\leq 1} \cong \wt{X}'_{\leq 1} \cong \Bl\left(X_{\leq
1},X_1\right).
$$
\item The resolutions $\pi:\wt{X} \to X$, $\pi':\wt{X}' \to X$ are
semismall with respect to the rank stratification on $X$.
\end{enumerate}
\end{remark}

\section{The $\C^*$-action}

We now turn to the proof of our main result, Theorem~\ref{main}. The
proof goes by a rather standard argument and uses a certain
$\C^*$-action on the quotient variety $X = V/G$ which has been very
well studied in many particular cases (see, e.g., \cite{quiver}).

We begin with some generalities. Let $Y$ be a smooth algebraic
variety equiped with an algebraic action of the group $\C^*$. For
every point $z \in Y$ denote by $\C^* \cdot z \subset Y$ the
$\C^*$-orbit of the point $x$, and let $\overline{\C^* \cdot z}
\subset Y$ be its Zariski closure. Choose a point $y \in Y$ fixed
under the $\C^*$-action, so that $\C^* \cdot y = \{y\}$. Consider
the the subset $S^0_+(y) \subset Y$ defined by
$$
S^0_+(y) = \left\{z \in Y\mid y \in \overline{\C^* \cdot
z}\right\},
$$
and let 
$$
S_+(y) = \overline{S^0_+(y)} \subset Y
$$
be it Zariski closure in $Y$. The subvariety $S_+(y) \subset Y$ is
called the {\em attraction subvariety} of the point $y \in Y$.

More generally, for every closed subvariety $Y_0 \subset Y$
consisting of points fixed by $\C^*$, define the attraction
subvariety $S_+(Y_0) \subset Y$ by
$$
S_+(Y_0) = \overline{\left\{z \in Y\mid \overline{\C^* \cdot z} \cap
Y_0 \neq \emptyset\right\}}.
$$

Let $Y_0 \subset Y$ be a connected component of the subvariety
$Y_{\C^*} \subset Y$ of points in $Y$ fixed by $\C^*$. Choose an
arbitrary point $y \in Y_0$.  Let $T_yY$ be the tangent space to $Y$
at the point $y \in Y$, and let
\begin{equation}\label{weights}
T_yY = \bigoplus_{p \in \Z} T^p_yY
\end{equation}
be its weight decomposition with respect to the $\C^*$-action: an
element $\lambda \in \C^*$ acts on $T^p_y$ by multiplication by
$\lambda^p$. Recall the following fact.

\begin{lemma}\label{standard}
\begin{enumerate}
\item The component $Y_0 \subset Y$ is smooth at $y \subset Y_0
\subset Y$, and the tangent space $Y_yY_0$ equals
$$
T_yY_0 = T^0_yY \subset T_yY = \bigoplus_{p \in \Z} T^p_yY.
$$
\item The attraction subvariety $S_+(Y_0) \subset Y$ is smooth at $y
\subset S_+(Y_0) \subset Y$, and the tangent subspace $T_yS_+(Y_0)
\subset T_yY$ equals
$$
T_yS_+(Y_0) = \bigoplus_{p \geq 0}T^p_yY \subset T_yY = \bigoplus_{p
\in \Z} T^p_yY.
$$
\end{enumerate}
\end{lemma}

\proof{} Consider the formal completion $\wh{Y}$ of $Y$ at the point
$y \in Y$. The group $\C^*$ does not act on the completion $\wh{Y}$,
but the differential of the action is a well-defined vector field
$\xi$ on $\wh{Y}$. The completion $\wh{Y}_0$ of the component $Y_0
\subset Y$ at $y \subset Y_0$ is the zero set of the vector field
$\xi$.  Moreover, the completion $\wh{S_+(Y_0)}$ of the attraction
variety $S_+(Y_0)$ at $y \in S_+(y)$ is a closed subvariety in
$\wh{Y}$, and it is defined by the ideal $I \subset \calo_{\wh{Y}}$
generated by formal function $f \in \calo_{\wh{Y}}$ with $\xi(f) =
0$.

Since the group $\C^*$ is reductive, the completion $\wh{Y}$,
equipped with the vector field $\xi$, is isomorphic to the
completion at $0$ of the tangent space $T_yY$, equipped with the
vector field defined by the natural $\C^*$-action on $T_yY$. Since
both completions $\wh{Y}_0$, $\wh{S_+(Y_0)}$ depend only on the
vector field $\xi$, it suffices to prove the lemma for $T_yY$
instead of $Y$. In this setting it is obvious.  \endproof

We also note the following obvious functorial property of the
attraction varieties.

\begin{lemma}\label{obvious}
If $f:Y \to Z$ is a $\C^*$-equivariant proper algebraic map between
two algebraic varieties equipped with $\C^*$-actions, then for any
closed subvariety $Y_0 \subset Y$ consisting of points fixed by
$\C^*$ we have
$$
f(S_+(Y_0)) \subset S_+(f(Y_0)). \qquad\qquad\qquad \square
$$
\end{lemma}

Return now to the situation of Theorem~\ref{main}, where $X = V/G$
is the quotient variety of the symplectic vector space
$$
V = V_o \oplus V_o^*
$$
by the finite group $G$ acting on $V_o$, and assume given a smooth
projective crepant resolution $\pi:\wt{X} \to X$. Let the group
$\C^*$ act on $V = V_o \oplus V_o^*$ by 
$$
\lambda \cdot \langle v, v' \rangle = \langle \lambda v, v' \rangle,
\qquad \lambda \in \C^*, \quad \langle v,v' \rangle \in V =
V_o \oplus V_o^*, \quad v \in V_o, v' \in V_o^*.
$$
This action commutes with the action of the finite group
$G$. Consequently, we obtain a $\C^*$-action on the quotient variety
$X = V/G$, which by Theorem~\ref{equivar}~\thetag{ii} lifts to a
$\C^*$-action on the resolution $\wt{X}$. We will call it the {\em
standard $\C^*$-action} on a crepant resolution $\pi:\wt{X} \to X$.

Note that the symplectic form $\Omega_V$ on $V = V_o \oplus V_o^*$
satisfies
\begin{equation}\label{weight.1}
\lambda^*(\Omega) = \lambda\Omega
\end{equation}
for every $\lambda \in \C^*$. By construction the same holds for the
symplectic form $\Omega$ on the open non-singular part $X_0 \subset
X$. The form $\Omega$ extends to a non-degenerate holomorphic form
$\Omega$ on the resolution $\wt{X}$ by
Proposition~\ref{crep=>sympl}, and the extension also satisfies
\eqref{weight.1} for every $\lambda \in \C^*$.

The crucial part of the proof of Theorem~\ref{main} is the following
fact.

\begin{prop}\label{discrete}
For every point $x \in X$ fixed by the $\C^*$-action, there exist at
most a finite number of points $\wt{x} \in \wt{X}$ fixed by the
$\C^*$-action and such that $\pi(\wt{x}) = x$.
\end{prop}

\proof{} By Proposition~\ref{tubular.res} and Theorem~\ref{formal}, it
suffices to prove the claim for the point $0 \in X$. Consider the
subvariety $\pi^{-1}(0)_{\C^*} \subset \pi^{-1}(0)$ of points in
$\pi^{-1}(0) \subset \wt{X}$ fixed by $\C^*$. Since the map
$\pi:\wt{X} \to X$ is proper, the subvariety $\pi^{-1}(0) \subset
\wt{X}$ is also proper, and it suffices to prove that
$\dim\pi^{-1}(0)_{\C^*} = 0$.

Let $Y \subset \pi^{-1}(0)$ be an irreducible component of the
variety $\pi^{-1}(0)_{\C^*}$, and let $y \in Y$ be an arbitrary
point in the non-singular part of the variety $Y \subset
\pi^{-1}(0)$. Consider the weight decomposition \eqref{weights} at
the point $y \in \wt{X}$. Note that the form $\Omega$ induces a
non-degenerate symplectic form on the tangent space $T_y\wt{X}$, and
the equation \eqref{weight.1} implies
$$
\dim T^p_y\wt{X} = \dim T^{1-p}_y\wt{X}, \qquad p \in \Z.
$$
Consider the attraction subvariety $S_+(Y) \subset \wt{X}$ of the
subvariety $Y \subset \wt{X}$. By Lemma~\ref{standard} we have
$$
\dim Y = \dim T^0_y\wt{X}
$$
and
$$
\dim S_+(Y_0) \geq \sum_{p \geq 0}\dim T^p_y\wt{X} = \dim Y +
\sum_{p > 0} \dim T^p_y\wt{X} = \dim Y + n,
$$
where
\begin{align*}
2n &= \sum_{p > 0} 2\dim T^p_y\wt{X} = \sum_{p > 0} \left(\dim
T^p_y\wt{X} + \dim T^{1-p}_y\wt{X}\right) \\
&= \sum_p \dim T_y^p\wt{X} = \dim\wt{X}.
\end{align*}
Now, by Lemma~\ref{obvious}, we have $\pi(S_+(Y)) \subset S_+(0)$,
and by the definition of the $\C^*$-action on $V = V_o \oplus V_o^*$
the attraction subvariety $S_+(0) \subset X$ of $0 \in X$ coincides
with the quotient
$$
S_+(0) = V_o/G \subset X = V/G
$$
of the subspace $V_o \subset V$ by the group $G$. Note that for
every subgroup $G_v \subset G$ the intersection $V_o \cap V^{G_v}$
coincides with the invariant vectors subspace $V_o^{G_v} \subset
V_o$. Therefore for every $k$ the intersection $S_+(0) \cap X_k$ of
the attraction variety $S_+(0) = V_+/G$ with the subvariety $X_k
\subset X$ of point of rank $k$ has the dimension
$$
\dim S_+(0) \cap X_k = \dim V_o - k = n - k.
$$
By Proposition~\ref{semismall}, this implies that 
$$
\dim \pi^{-1}(S_+(0)) \leq n.
$$
Since $\pi(S_+(Y))$ lies within $S_+(0)$, we conclude that
$$
n \geq \dim \pi^{-1}(S_+(0)) \geq \dim S_+(Y) \geq \dim Y + n,
$$
which yields the required equality $\dim Y = 0$.
\endproof

We can now prove Theorem~\ref{main}. The argument more or less
literally repeats the proof of Theorem~4.2 in the paper
\cite{part-res}.

\proof[Proof of Theorem~{\normalfont\ref{main}}.] Consider the
subvariety $\wt{X}_{\C^*} \subset \wt{X}$ of the points fixed by the
standard $\C^*$-action on the crepant resolution $\pi:\wt{X} \to
X$. Since the variety $\wt{X}$ is smooth, the subvariety
$\wt{X}_{\C^*} \subset X$ is a union of smooth connected components.

By definition of the standard $\C^*$-action on $X$, the quotient
$X_o = V_o^* / G \subset X$ of the $G$-invariant Lagrangian subspace
$V_o^*$ by the group $G$ is the subvariety of $\C^*$-fixed point in
the quotient $X = V/G$. The generic point of the subvariety $X_o
\subset X$ lies in the dense open subset $X_0 \subset X$ of vectors
of rank $0$. Since the projection $\pi:\wt{X} \to X$ is one-to-one
over the subset $X_0 \subset X$, there exists a connected component
$Y \subset \wt{X}_{\C^*}$ of the fixed-points subvariety
$\wt{X}_{\C^*} \subset \wt{X}$ such that the induced map $\pi:Y \to
X_o$ is dominant and generically one-to-one.

By Proposition~\ref{discrete}, the dominant map $\pi:Y \to X_o$ is
in fact finite. But the quotient variety $X_o = V_o^* / G$ is
normal. Therefore the finite map $\pi:Y \to X_o = V_o^* / G$ is not
only generically one-to-one, but induces an isomorphism between $Y$
and $X_o$. This implies that the quotient $X_o = V_o^* / G$ is
smooth.

To finish the proof of Theorem~\ref{main}, it remains to invoke the
classic Theorem~\ref{bourb}.
\endproof

\end{document}